\documentclass[12pt]{amsart}
\usepackage{amsmath} 
\usepackage{amsfonts,amssymb}
\usepackage{array,delarray}
\usepackage[dvips]{graphicx}

\newcommand{\mm}{\mathcal{M}}
\newcommand{\kk}{\mathbb{K}}
\newcommand{\acyc}{\mathrm{Acyc}}

\newtheorem{theorem}{Theorem}[section]   
\newtheorem{cor}[theorem]{Corollary}     
\newtheorem{lemma}[theorem]{Lemma}         
\newtheorem{proposition}[theorem]{Proposition}  

\theoremstyle{definition}
\newtheorem{definition}[theorem]{Definition}   

\theoremstyle{remark}
\newtheorem{remark}[theorem]{Remark}        
\newtheorem{example}[theorem]{Example}        

\numberwithin{equation}{section}     

\begin{document}

\title[Finite Dynamical Systems]  	
{Classification of Finite Dynamical Systems} 

\author[Garcia et. al.]{Luis Garcia  \and Abdul Salam Jarrah 
\and Reinhard Laubenbacher}

\address[Luis Garcia]
{Department of Mathematical Sciences\\
New Mexico State University\\
Las Cruces, NM 88003}    
\email{galuis@nmsu.edu}

\address[Abdul Salam Jarrah]
{Department of Mathematical Sciences\\
New Mexico State University\\
Las Cruces, NM 88003}    
\email{ajarrah@nmsu.edu} 

 \thanks{
This work was partially supported by funds from a partnership initiative
between NMSU and Los Alamos National Laboratory. 
}

\address[Reinhard Laubenbacher]
{Department of Mathematical Sciences\\
New Mexico State University\\
Las Cruces, NM 88003}
\email{reinhard@nmsu.edu}

   
\begin{abstract}
This paper is motivated by the theory of sequential dynamical systems,
developed as a basis for a mathematical theory of computer simulation.
It contains a classification of finite dynamical systems on 
binary strings, which are obtained by composing functions defined on the
coordinates. The classification is in terms of the dependency relations
among the coordinate functions.  It suggests a natural
notion of the linearization of a system.  Furthermore, it contains a sharp
upper bound on the number of systems in terms of the dependencies among
the coordinate functions.  This upper bound generalizes an upper bound
for sequential dynamical systems.
\end{abstract}

\maketitle

\section{Introduction}

In this paper a {\it finite dynamical system} will mean a mapping
from the set of binary strings of a certain length to itself.
Such mappings occur in a variety of contexts, in particular in
the theoretical study of computer simulations.  Representing
a computer simulation as a particular finite dynamical system is one possible
approach to its mathematical analysis,
carried out in \cite{br, bmr,bmr2}.  
The finite dynamical systems considered there, so-called {\it sequential
dynamical systems} (SDS), incorporate the essential features of a computer simulation.
Local variables $v_1,\ldots ,v_n$ 
take on binary states which evolve in discrete
time, based on a local update function $f^i$ attached to each variable $v_i$, and which depends
on the states of certain other variables, encoded by the edges of a 
{\it dependency graph} $Y$ on the vertices $v_1,\ldots ,v_n$.  
Finally, an update schedule prescribes how
these local update functions are to be composed in order to generate a global update
function
$$
f:\{0,1\}^n\longrightarrow \{0,1\}^n
$$
of the system.  
An important question, which can be answered in this setting, is how many different
systems one can generate simply by varying the update schedule.  The upper bound is
in terms of invariants of the dependency graph $Y$.

In applications the dependency graph $Y$ frequently varies over time, however.
The need for a framework that allows for such a change inspired the investigation
of properties of tuples of ``local'' functions in \cite{LP} and certain equivalence
relations on them. That paper also contains a Galois correspondence between sets of
tuples of local functions and certain graphs.  Tuples of local functions can
be interpreted as parallel systems $f:\{0,1\}^n\longrightarrow \{0,1\}^n$, so that
the results pertain to the study of parallel systems as well.
The present paper makes the
connection between tuples of local functions (parallel systems)
and sequential systems, in particular SDS,
explicit by exploiting this Galois correspondence.  In order to
describe our main results we need to recall some definitions and results from \cite{LP}.

\medskip 
Let $\kk=\{0,1\}$, 
and let $\kk^n$ be the $n$-fold cartesian product of $\kk$.  
When convenient we will view $\kk$ as the field with two elements.

\begin{definition}
Let $n$ be a positive integer, let $d$ be a nonnegative integer, and
let $Y$ be a graph with vertex set $\{1,\dots,n\}$.
\begin{enumerate}
\item A function $f: \kk^{n} \longrightarrow \kk^{n}$ is \emph{$d$-local} on Y if,
for any $1 \leq j \leq n$, the $j$-th coordinate of the value of $f$
on $x \in \kk^{n}$ depends only on the value of those coordinates of
$x$ that have distance less than or equal to $d$ from vertex $j \in
Y$. In other words, if $f(x) = (f_{1}(x),\dots,f_{n}(x))$, then
$f_{j}: \kk^{n} \longrightarrow \kk$ depends only on those coordinates that have
distance less than or equal to $d$ from $j$.

\item For $1\leq d < n$ and $1\leq j\leq n$, let $L_{d}^{j}(Y)$ be the set of 
all functions $f: \kk^{n} \longrightarrow \kk^{n}$ such that 
\begin{equation*}
f(x_{1},\dots,x_{n}) = (x_{1},\dots,x_{j-1},f_{j}(x),x_{j+1},\dots,x_{n}),
\end{equation*}
and $f_{j}:\kk^{n}\to \kk$ depends only on the values of those
coordinates of $x$ which have distance at most $d$ from $j$ in
$Y$. Hence $L_{d}^{j}(Y)$ consists of $d$-local functions on
$\kk^{n}$, which are the identity on all but possibly the $j$-th
coordinate.

\item For $d=n$, define $L_{n}^{j}(Y)$ to be the set of all functions on $\kk^{n}$,
which are the identity on all but possibly the $j$-th coordinate. Observe
that if $Y$ is connected, then this definition of $L_{n}^{j}(Y)$
directly extends the definition in $(2)$. 
\end{enumerate}
\end{definition}

Observe that  $L_{0}^{j}(Y)=L_0^j$ does not depend on the graph $Y$ and 
neither does $L_{n}^{j}(Y)$. 
Furthermore, $L_{0}^{j}$ is isomorphic to ${\rm Map}(\kk,\kk)$,
that is, it contains all four possible functions, namely the identity
on $\kk$, the two projections to one element in $\kk$, and the
inversion.

In this paper we study the set 
\[ L_{n}^{1}\times\dotsm\times L_{n}^{n} = \{(f^{1},\dots,f^{n})\mid
f^{i} \in L_{n}^{i}\},\]
that is, the set of $n$-tuples of functions $f^i:\kk^n\longrightarrow \kk^n$
which only change the $i$th coordinate.  To be precise, $f^{i}(x) = (x_{1},\dots,x_{i-1},f_{i}^{i}(x),x_{i+1},\dots,x_{n})$, with
arbitrary functions $f_{i}^{i}:\kk^{n}\longrightarrow \kk$. We denote by
$\mathcal{F}$ the power set of this set without the empty set. The
following theorem is one of the main results in \cite{LP}.

\begin{theorem}\label{thm:galois}
There is a Galois correspondence between $\mathcal{F}$ and the set
$\mathcal{G}$ of subgraphs of the complete graph $K_{n}$ on the vertex
set $\{1,\dots,n\}$.
\end{theorem}

For the convenience of the reader we 
recall the construction of this Galois correspondence. Let $F \in
\mathcal{F}$. Define a subgraph $\Phi(F)$ of $K_{n}$ as follows. First
construct the set $\widetilde{F}$ of all $n$-tuples $\widetilde{f} =
(\widetilde{f^{1}},\dots,\widetilde{f^{n}})$, which either are in $F$ or arise
from an element in $F$ by replacing one of the coordinates by a
$0$-local function, that is, by a function from $L_{0}^{i}$ for some
$i$. Now define the graph $\Phi(F)$ as follows. An edge $(i,j)$ of
$K_{n}$ is in $\Phi(F)$ if and only if
$\widetilde{f^{i}}\circ\widetilde{f^{j}} = \widetilde{f^{j}}\circ\widetilde{f^{i}}$ 
for all 
$\widetilde{f} = (\widetilde{f^{1}},\dots,\widetilde{f^{n}}) \in \widetilde{F}$. 

Conversely, let $G \subset K_{n}$ be a subgraph. We define a set
$\Psi(G)$ of $n$-tuples of functions on $\kk^{n}$ by 
\[
\Psi(G) = L_{1}^{1}(\overline{G}) \times  L_{1}^{2}(\overline{G})\times \dotsm
 \times L_{1}^{n}(\overline{G}),
\]  
where $\overline{G}$ is the complement of $G$ in $K_{n}$.
Then $\Phi$ and $\Psi$ together form the desired Galois correspondence.

In particular, if $F\in\mathcal F$ consists of one element $f=(f^1,\ldots ,f^n)$,
then the graph $\Phi(\{f\})=\Phi(f)$ encodes the dependency relations among the 
local functions $f^i$.  Conversely, for a subgraph $G\subset K_n$, the set
$\Psi(G)$ contains all $n$-tuples of local functions whose dependency relations
are modeled by $G$.  This observation provides the paradigm for the
results in this paper.

First we encode the local functions as polynomials, which allows us to give
an algebraic criterion to compute $\Phi(f)$.  More importantly, it suggests
a natural choice for the linearization of a system, which we define in the
next section.  Using this notion of linearization we explore a natural 
equivalence relation on tuples, setting two equivalent if they have the same
dependency graph.  We show that systems that are equivalent in this sense
have the same linearization.

Finally, we consider systems $f:\kk^n\longrightarrow \kk^n$ which are obtained
by composing local functions from an $n$-tuple $f=(f^1,\ldots ,f^n)$.  We
show that, if $t\geq 1$ is an integer and $W_t$ is the set of all words in the
integers $1,\ldots ,n$ of length $t$, allowing for repetitions and for the case
that $t<n$, then we obtain an upper bound for the number of different
systems $f^{\pi}:\kk^n\longrightarrow \kk^n$ one can construct by forming
$$
f^{\pi}=f^{i_t}\circ\cdots\circ f^{i_1},
$$
where $\pi =(i_1,\ldots ,i_t)$ ranges over all elements in $W_t$.  

This upper bound generalizes one for SDS, derived in \cite{R}.  It suggests that
a part of the theory of SDS can be derived for systems that are SDS-``like,'' but
have fewer restrictions on the local functions.  In particular, it is not 
necessary to make the dependency graph an explicit part of the data defining an
SDS.  This approach is explored further in \cite{LP2}, in which a
more general notion
of SDS is introduced, and morphisms of SDS are defined, forming
a category with interesting properties, that contains ``classical'' SDS as a 
subcategory.   A morphism between two SDS can be viewed as a simulation of
one system by the other.

\section{Computation of $\Phi$}

In this section we give a method for computing the graph $\Phi(f)$
for an $n$-tuple of local functions $f$.  It relies on the 
representation of local functions as polynomials, for which we now
give an elementary proof.

\begin{lemma}
Let $f:\kk^n\longrightarrow \kk$ be a function.  Then $f$ can be
represented as a polynomial.  That is, there is a polynomial
$p\in \kk[x_1,\ldots ,x_n]$ such that
$$
f(a_1,\ldots ,a_n)=p(a_1,\ldots ,a_n),
$$
for all $(a_1,\ldots ,a_n)\in\kk^n$.
\end{lemma}

\begin{proof}
We will prove the lemma by showing that there are exactly as many 
different polynomial functions as there are functions.  Since
$|\kk^n|=2^n$, there are $2^{2^n}$ functions from $\kk^n$ to $\kk$.
Now let $V\subset\kk[x_1,\ldots ,x_n]$ be the subspace with the basis
containing the $2^n$ monomials
$$
\mathcal B=\{1,m_1,m_2,\ldots ,m_{2^n-1}\},
$$
where $m_i=x_1^{b_1}\cdots x_n^{b_n}$ is given 
by the binary expansion $b_1\ldots b_n$ of $i$.
Then $V$ contains $2^{2^n}$ elements, which is equal to the number of functions
$\kk^n\longrightarrow \kk$.  Thus it suffices to show that different polynomials
give rise to different functions on $\kk^n$.

Suppose that $f,g\in V$, with $f\neq g$.  
Let $m=x_1^{a_1}\cdots x_n^{a_n}$ be a monomial of smallest total degree which
is in $f$ but not in $g$, and let $h=f-g=m+q$.  Let ${\bf a}=(a_1,\ldots ,a_n)$.  
Then $h({\bf a})=1+q({\bf a})$.  Now observe that any monomial of $q$ which involves a
variable that does not appear in $m$ evaluates to $0$ at ${\bf a}$.  But since $m$ was
chosen to have minimal degree, there cannot be any monomials in $q$ that involve only
variables appearing in $m$.  Hence $q({\bf a})=0$, and $h({\bf a})=1$.  This shows that
$f\neq g:\kk^n\longrightarrow \kk$, and the proof is complete. 
\end{proof}

Therefore, any function in $L_{n}^{j}$ can be written as an
$n$-tuple of polynomials in $n$
variables $x_{1},\dots, x_{n}$. Also observe that $(x_j)^m = x_j$ for all 
$1\leq j\leq n$ and for all $m>0$. Furthermore, 
$\overline{x_j}=1 + x_j$. 
Thus, we can represent an element of $L_{n}^{1} \times\dotsm\times L_{n}^{n}$ 
as an $n\times (2^n)$-matrix in which the $i$th row corresponds to the $i$th
local function $f_{i}^i$. 

\begin{example}
Let $f = (f^1, f^2, f^3)$ where
\begin{align*}
f^1 &= (1 + x_3 + x_1x_2,x_{2},x_{3}), \\
f^2 &= (x_{1},x_2 + x_1x_3 + x_2x_3,x_{3}), \\
f^3 &= (x_{1},x_{2},1 + x_1x_2x_3).
\end{align*}

Then the $3\times 8$ matrix, denoted by $M_{f}$, associated to $f$ with respect to the
ordered basis 
\[
\{1, x_1, x_2, x_3, x_1x_2, x_1x_3, x_2x_3, x_1x_2x_3\}, 
\]
 is equal to 
\[
\begin{array}({cccccccc}) 
1 & 0 & 0 & 1 & 1 & 0 & 0 & 0\\
0 & 0 & 1 & 0 & 0 & 1 & 1 & 0\\
1 & 0 & 0 & 0 & 0 & 0 & 0 & 1 
\end{array}.
\]
\end{example}

Let $f = (f^1, \dots,f^n):\kk^n \longrightarrow \kk^n$
 be an element of $L_{n}^{1}\times\dotsm\times L_{n}^{n}$  
and let $\Phi(f)$ be its 
associated graph. The following result gives an important characterization
of $\Phi(f)$ in terms of the polynomial representation of the entries
of $f$.

\begin{proposition}\label{main:thm1}
There is an edge between vertex $i$ and vertex $j$ in  
$\Phi(f)$ if and only if $x_i$ does not divide 
any monomial of $f^j_j$ and $x_j$ does not divide any monomial of 
$f^i_i$.
\end{proposition}

\begin{proof}
Suppose first that $x_i$ does not appear in $f_j^j$ and $x_j$ does not
appear in $f_i^i$.  It is clear that then $f^i$ and $f^j$ commute.
Now consider an $\widetilde{f}$ which is obtained from
$f$ by replacing the $j$th coordinate by a $0$-local function $\widetilde f^j$.
It also does not depend on $x_i$.
Similarly, no $\widetilde{f^i}$ depends on $x_j$. 
Thus 
$\widetilde{f^j} \circ \widetilde{f^i} = \widetilde{f^i} \circ
\widetilde{f^j}$ for all $\widetilde{f} =
(\widetilde{f^1},\dots,\widetilde{f^n})$ obtained from $f$ by replacing
one of the coordinates by a $0$-local function.
Hence the edge $(i,j)$ is in the graph $\Phi(f)$.

Conversely, suppose that $x_i$ divides 
a monomial of $f^j_j$ or $x_j$ divides a monomial of $f^i_i$. Without 
loss of generality, suppose  $x_i$ divides a monomial of $f^j_j$.
Let $\mm = \{ m_1, \dots, m_t \}$ be the set of all monomials of $f^j_j$ 
such that $x_i | m_l$ for all $l=1, \dots, t$. Let $m_s$ be a
monomial in $\mm$ of minimal degree,
say $m_s = x_{s_1}x_{s_2}\dotsm x_{s_r}$, where
$x_i \in \{ x_{s_1}, \dots,x_{s_r}\} \subseteq \{x_1,\dots,x_n\}$. Define 
$a=(a_1,\dots, a_n) \in \kk^n$ as follows, for $1 \leq l \leq n$,
\begin{equation*}
a_l = \left \{ \begin{array}{ll}
		1 & \mbox{if $l \in \{s_1, \dots, s_r \}$,} \\
		0 & \mbox{otherwise.}
		\end{array}
	\right.
\end{equation*}
Then the $j$th coordinate of $f^j(a)$ is equal to
$$
f^j_j(a)= (m_1 + \cdots + m_t)(a) = m_s(a) = 1.
$$
Let $\widetilde{f^i}$ be the $0$-local function, projection to zero, in the $i$th
coordinate. The $j$th coordinate of $\widetilde{f^i} \circ f^j(a)$ 
is equal to
$$
(\widetilde{f^i} \circ f^j(a))_j = f^j_j(a) = 1.
$$

On the other hand, 
$\widetilde{f^i}(a) = (a_1, \dots,a_{i-1},0,a_{i+1},\dots,a_n)$. 
So the $j$th coordinate of $f^j \circ \widetilde{f^i}(a)$ is  

\begin{eqnarray*}
(f^j \circ \widetilde{f^i}(a))_j &=& f^j_j(a_1, \dots,a_{i-1},0,a_{i+1},\dots,a_n)\\
	&=& (m_1 + \cdots + m_t)(a_1, \dots,a_{i-1},0,a_{i+1},\dots,a_n)\\
	&=& 0. 
\end{eqnarray*}

Thus 
$(\widetilde{f^i} \circ f^j(a))_j \neq (f^j \circ \widetilde{f^i}(a))_j$.
Therefore, $\widetilde{f^i} \circ f^j \neq f^j \circ \widetilde{f^i}$, and
hence there cannot be an edge between vertices $i,j$ in $\Phi(f)$.  
\end{proof}

This proposition provides an easy algorithm to compute $\Phi(f)$ for 
any $f$, or, more generally, any set of such functions.  (A C++ implementation
is available from the authors.)  It also suggests a natural notion
of a linear system, which will be explored in the next section.
\section{Linearization of Systems}

For classical dynamical systems, a standard technique is to linearize the
system, and then study the linearization.  Proposition \ref{main:thm1}
suggests a natural definition of linearity for a finite system
$$
f=(f^1,\ldots ,f^n):\kk^n\longrightarrow \kk^n.
$$
One can then define the notion of linearization for a finite system
and study the relationship between the original system and its 
linearization.  In this section we study the notion of linearization
with respect to the dependency relations among the local functions.

Each $f^i:\kk^n\longrightarrow \kk^n$ can be represented as a polynomial in
the variables $x_1,\ldots ,x_n$, which changes only the $i$th coordinate.
That is, $f^i=(f^i_1,\ldots ,f^i_n)$, with each $f^i_j\in \kk[x_1,\ldots ,x_n]$.

\begin{definition}
A system $f = (f^1,\dots,f^n) :\kk^n \longrightarrow \kk^n$ is called 
{\em linear} if all functions $f_{i}^{i}$ are $\kk$-linear polynomials.
\end{definition}

Observe that a linear system $f=(f^{1},\dots,f^{n})$ can be
represented by an $(n\times n)$-matrix with entries in $\kk$, 
since the constant term and all
nonlinear terms of all functions $f_{j}^{i}$ are equal to zero.  Conversely,
any $(n\times n)$-matrix $M$ over $\kk$ can be interpreted as a linear system
$$
l_M:\kk^n\longrightarrow \kk^n .
$$

\begin{lemma}\label{main:cor1}
Let $G$ be a graph on $n$ vertices.  Then there exists a linear system 
$l_G = (l^1,\dots,l^n):\kk^n \longrightarrow \kk^n$ such that $\Phi(l_G) = G$.
\end{lemma}

\begin{proof}
We construct a linear system $l_G$ by providing the associated 
$n\times n$-matrix $M_{l_G}$ of $l_G$. For $1 \leq i < j \leq n$, let
\[
(M_{l_G})_{ij} = (M_{l_G})_{ji} = \left \{ \begin{array}{ll}
				1 & \mbox{if the edge $(i,j) \notin E(G)$}, \\
				0 & \mbox{otherwise.}
			\end{array} \right.
\]
And let $(M_{l_G})_{ii} = 0$, for all $1\leq i \leq n$. It is clear that 
$l_G$ is a linear system. Moreover, by Proposition \ref{main:thm1}, $ \Phi(l_G) = G$. 
\end{proof}

\begin{remark}\label{main:rem1}
The matrix  $M_{l_G}$ is the adjacency matrix of the complement of 
the graph $ \Phi(l_G) = G$.  Observe that $l_G$ is not the only linear system
which gives the graph $G$.  The matrix $M_{l_G}$ is symmetric.  Consider a pair
$i\neq j$ for which the $(i,j)$-entry (and the $(j,i)$-entry) is equal to $1$.  If
we change exactly one of these two entries to a $0$, then 
$l_G$ and $l_G'$ have the same graph $G$ as their image under $\Phi$.  We will see later
(Theorem \ref{equivalence}),
however, that the incidence matrix of the complement of $G$ is not just a 
canonical, but a natural choice of linear system in the inverse image
of a graph $G$ under $\Phi$.
\end{remark}

\begin{theorem}\label{main:thm2}
Let $f = (f^1,\dots,f^n) :\kk^n \longrightarrow \kk^n$ be a system and
let $G = \Phi(f)$.
There exists a linear system 
$l_G = (l^1,\dots,l^n) :\kk^n \longrightarrow \kk^n$ such that
$\Phi(l_G) = \Phi(f)$.
\end{theorem}

\begin{proof}
By Lemma \ref{main:cor1}, there exists a 
linear system $l_G$ such that $\Phi(l_G) = G = \Phi(f)$.
\end{proof}

\begin{definition}
The linear system $l_G$ corresponding to the adjacency matrix of the
complement of $G$ is called the {\em linearization}
of the system $f$.
\end{definition}

The linearization of a system should ideally bear
a certain relationship to the original system.  Our focus in this paper is
the structure of the dependencies among the component functions.  The
linearization of a system as defined here, has the same dependency structure
among the component functions.  Its dynamics will in general be very different, however.


\section{Graph equivalence}
Using the Galois correspondence we now define an
equivalence relation on $n$-tuples of local functions, which captures 
equivalence of dependency relations among the entries in the tuples.

\begin{definition}
We say that $(f^{1},\dots,f^{n}),(g^{1},\dots,g^{n}) \in
L_{n}^{1}\times \dotsm \times L_{n}^{n}$ are \emph{graph equivalent}
if and only if $\Phi((f^{1},\dots,f^{n}))$ is isomorphic to
$\Phi((g^{1},\dots,g^{n}))$. This relation is denoted by
$(f^{1},\dots,f^{n})\sim (g^{1},\dots,g^{n})$.
\end{definition}

\begin{example}
Let $f = (f^{1},f^{2},f^{3})$, where 
\begin{align*}
f^{1}(x_{1},x_{2},x_{3})&=(x_{2}+x_{3},x_{2},x_{3}),\\
f^{2}(x_{1},x_{2},x_{3})&=(x_{1},1+x_{2},x_{3}),\\ 
f^{3}(x_{1},x_{2},x_{3})&=(x_{1},x_{2},x_{1}).
\end{align*}

We will usually denote such a function by 
$f = (x_{2}+x_{3},1+x_{2},x_{1})$.  Let $g=(x_{2},x_{3},0)$. Then 
$f\sim g$, since their graphs both have three vertices $1,2,3$, 
with edge $(2,3)$ in $\Phi(f)$, and edge $(1,3)$ in $\Phi(g)$. 
\end{example}

\medskip

The main result in this section is a characterization of graph equivalence
of systems using a relationship on the matrices of the corresponding
linearizations.  We first define an action of the symmetric group on 
the set of matrices.

\begin{definition}\label{def:action}
Let $S_{n}$ be the group of permutations of $n$
elements. We define an $S_{n}$-action on the set of $n\times
n$-matrices as follows. For any $(n\times n)$-matrix $M$ and $\pi \in S_n$,
$\pi M$ is the  $(n\times n)$-matrix such that
\[
(\pi M)_{ij} = M_{\pi^{-1}(i)\pi^{-1}(j)}.
\]
That is, $\pi$ acts on $M$ by permuting rows {\it and} columns.
\end{definition}

\begin{proposition}
The $S_{n}$-action in Definition \ref{def:action} is a group action of $S_n$
on the set of $(n\times n)$-matrices over $\kk$.
\end{proposition}

\begin{proof}
Clearly $e M = M$, where $e$ is the identity permuation. Let $\pi,\sigma \in S_{n}$, then
\begin{align*}
(\pi(\sigma M))_{ij} &= (\sigma M)_{\pi^{-1}(i)\pi^{-1}(j)} =
M_{\sigma^{-1}\pi^{-1}(i)\sigma^{-1}\pi^{-1}(j)}\\
&= M_{(\pi\sigma)^{-1}(i)(\pi\sigma)^{-1}(j)} = (\pi\sigma)M_{ij}. 
\end{align*}
\end{proof}

A graph automorphism can be represented by a permutation
of the vertices of the graph which preserves adjacency. 
So we can represent automorphisms of a
graph with $n$ vertices as permutations in $S_{n}$, and obtain in this
way an $S_n$-action on the set of subgraphs of $K_n$. 

\begin{theorem}\label{main:thm3}
Let $G_1$ and $G_2$ be two graphs with $(n\times n)$-adjacency matrices $M_1$ 
and $M_2$, 
respectively. Then $G_1$ is isomorphic to $G_2$ via 
an automorphism $\pi \in S_{n}$ 
(i.e., $\pi(G_1) = G_2$) if and only if $\pi M_1 =M_2$.
\end{theorem}

\begin{proof}
Suppose first that $\pi M_1 =M_2$.
For any two vertices $i$ and $j$ of $G_2$, 
$(i,j) \in E(G_2)$ if and only if  
$$
(M_2)_{ij} =(M_2)_{ji} =1.
$$
This happens if and only if
$$
(M_1)_{\pi^{-1}(i)\pi^{-1}(j)} = (M_1)_{\pi^{-1}(j)\pi^{-1}(i)} =1,
$$
which is the case if 
and only if the edge $(\pi^{-1}(i),\pi^{-1}(j)) \in E(G_1)$. Thus $\pi(G_1) = G_2$. 

Conversely, if $\pi(G_1) = G_2$, then 
$(M_2)_{ij} = 1$ if and only if $(i,j) \in E(G_2)$, that is,
if and only if $(\pi^{-1}(i),\pi^{-1}(j)) \in E(G_1)$.  This is the case
 if and only if
$(M_1)_{\pi^{-1}(i)\pi^{-1}(j)} = (M_1)_{\pi^{-1}(j)\pi^{-1}(i)} =1$.
Therefore, $\pi \cdot M_1 =M_2$.   
\end{proof}

\begin{theorem}\label{equivalence}
Let $f$ and $g$ be two systems on $\kk^n$. Then $f$ is graph 
equivalent to $g$ (i.e., $\pi(\Phi(f)) = \Phi(g)$ for some
$\pi \in S_n$) if and only if $\pi \cdot M_{l_{\Phi(f)}} =  M_{l_{\Phi(g)}}$.
\end{theorem}

\begin{proof}
First observe that if $G_1$ and $G_2$ are graphs, then
$G_1$ is isomorphic to $G_2$ via $\pi$ if and only if  
$\overline{G_1}$ is isomorphic to  $\overline{G_2}$ via $\pi$.
Therefore, $\pi(\Phi(f)) = \Phi(g)$ if and only if  
$\pi(\overline{\Phi(f)}) = \overline{\Phi(g)}$. By Remark \ref{main:rem1},
$M_{l_{\Phi(f)}}$ (resp. $M_{l_{\Phi(g)}}$) is the adjacency 
matrix of $\overline{\Phi(f)}$
(resp. $\overline{\Phi(g)}$). 
Now, by Theorem \ref{main:thm3},  
$\pi(\overline{\Phi(f)}) = \overline{\Phi(g)}$ if and only if
$\pi \cdot M_{l_{\Phi(f)}}= M_{l_{\Phi(g)}}$.
Therefore,  $\pi(\Phi(f)) = \Phi(g)$ if and only if 
$\pi \cdot M_{l_{\Phi(f)}} =  M_{l_{\Phi(g)}}$.
\end{proof}

This theorem suggests that our choice of the adjacency matrix of
$\overline{\Phi(f)}$ as the linearization of a system $f$ is a natural one,
since it preserves graph equivalence and makes the construction of $l_G$ from
$G$ equivariant with respect to the $S_n$-action on both.
\section{An Upper Bound for Sequential Systems}

The real object of interest in many cases are composed systems $f:\kk^n\longrightarrow \kk^n$
rather than merely tuples of local functions.  Such systems are obtained by
composing the local functions in some order.  The order often corresponds to a choice
of update schedule of the variables in a system, such as a simulation.  A
theoretical question which has important practical consequences is how many different
systems one can obtain by simply varying the update schedule of the variables, that is,
by composing the local functions in a different order.  In this section we derive
an upper bound for this number.  

\begin{definition}
Let $f = (f^1,\dots,f^n) : \kk^n \to \kk^n$, and let $W_t$ be
the set of all words on $\{1,\dots,n\}$ of length $t$, for some $t\geq 1$, allowing 
for repetitions. For 
$\pi =(i_1,\ldots ,i_t)\in
W_t$, we denote by $f^{\pi}$ the finite dynamical system given by 
\[ 
f^{i_t}\circ \cdots \circ f^{i_1}: \kk^n\to\kk^n.
\]
 
Let $F_{W_t}(f) = \{ f^{\pi} \mid \pi \in W_t \}$, the collection of all systems
$\kk^n\longrightarrow \kk^n$ that can be obtained by composing the coordinate functions
of $f$ in all possible ways, using up to $t$ of them.
\end{definition}

We now define an equivalence relation on $W_t$. 

\begin{definition}
Let $G$ be a graph on the $n$ vertices $1,\ldots ,n$. 
Let $\sim_G$ be the equivalence relation
on $W_t$ generated by the following relation.  Let $\pi=(i_1,\ldots ,i_t)\in W_t$. 
For $1\leq k<t$, if $i_k=i_{k+1}$ or there is no 
edge between $i_k$ and $i_{k+1}$ in $G$, then
$$
\pi\sim_G\pi',
$$
where
$
\pi'=(i_1,\ldots ,i_{k+1},i_k,i_{k+2},\ldots ,i_t).
$
\end{definition}

\begin{remark}
If $f = (f^1,\dots,f^n) : \kk^n \to \kk^n$ is  
such that $\Phi(f) = G$, and $\pi \sim_{\bar{G}} \pi'$, then 
\[f^{i_t}\circ \cdots \circ f^{i_k}\circ f^{i_{k+1}} \circ \cdots \circ f^{1} = 
f^{i_t}\circ \cdots \circ f^{i_{k+1}}\circ f^{i_k} \circ \cdots \circ f^{1}.
\]
\end{remark}




We now derive an upper bound on the size of the set $F_{W_t}(f)$, that is,
on the number of different systems one obtains by composing the coordinate
functions of $f$ in all possible orders, with up to $t$ of them at a time.
 
\begin{definition}
Let $f=(f^1,\ldots ,f^n):\kk^n\longrightarrow \kk^n$, and let 
$G(f)=G=\Phi(f)$.
Let $\pi =(i_1,\ldots ,i_t)\in W_t$. Let $H_{\pi}(f)=H_{\pi}$ be the
graph on $t$ vertices $v_1,\ldots ,v_t$, corresponding to 
$i_1,\ldots ,i_t$ (with $v_a\neq v_b$ even in the case that $i_a=i_b$), 
with an edge between $v_a$ and $v_b$ 
if and only if the following two conditions
hold:
\begin{enumerate}
\item $i_a \neq i_b$,
\item the edge $(i_a,i_b)$ is not in $G$.
\end{enumerate}
\end{definition}

\begin{remark}
Observe that if  $\pi \in S_n$, that is, $t=n$ and $\pi$ contains
no repetitions, then $H_{\pi} = \bar{G}$.
\end{remark}

Let $\acyc(H)$ be the set of all acyclic orientations of a graph
$H$.  Given $\pi=(i_1,\ldots ,i_t)\in W_t$, we construct an acyclic
orientation of $H_{\pi}$ by orienting an edge $(v_i,v_j)$ toward the
vertex whose label occurs first in $\pi$.  If all entries of $\pi$ are
distinct, then this clearly produces an acyclic orientation.  But
even if an entry is repeated we cannot produce an oriented cycle,
since there is no edge between the vertices corresponding to the
repetitions.  Denote this acyclic orientation by 
$\mathcal O_{\pi}(f)$.

\begin{lemma}\label{orientation}
If $\pi\sim_{\overline{G}}\pi'$, 
then $H_{\pi}(f)=H_{\pi'}(f)$ and $\mathcal O_{\pi}(f)=\mathcal O_{\pi'}(f)$.
\end{lemma}

\begin{proof}
If $\pi\sim_{\overline{G}}\pi'$, then they differ by a sequence of transpositions of adjacent
letters, which are either equal, or for which the corresponding 
vertices in $G$ are connected by an edge.
Hence $H_{\pi}(f)$ and $H_{\pi'}(f)$ have the same vertex set.  
Furthermore, an edge $(a,b)$ is in $H_{\pi}(f)$ if and only
if $i_a\neq i_b$ and $(i_{a},i_{b})$ is an edge in $G$.  Similarly for $H_{\pi'}(f)$.
Observe that the transposition in $\pi$ of adjacent letters which are connected 
by an edge in $G$ does not change the resulting acyclic orientation,
because, by construction, the vertices $v_{a}$ and $v_{b}$ are not
connected by an edge in $H_{\pi}(f)$. Hence the proof of
the lemma is complete.
\end{proof}

The next proposition is a generalization of a result from Cartier-Foata normal form 
theory.  See, e.g., \cite{D,R}.

\begin{proposition}\label{main:lem1}
Let $f$ be a system and $G=\Phi(f)$.
There is a one-to-one correspondence
$$
\psi_G:W_t/\sim_{\overline{G}}\longrightarrow \{{\rm Acyc}(H_{\pi}(f))|\pi\in W_t\}.
$$
\end{proposition}

\begin{proof}
We assign to a word $\pi\in W_t$ the associated acyclic orientation 
$\mathcal O_{\pi}(f)$
on $H_{\pi}(f)$.  
By Lemma \ref{orientation} this induces a mapping 
$\psi_G$ on $W_t/\sim_{\overline{G}}$.
There is an obvious inverse mapping, assigning to an acyclic orientation on $H_{\pi}$
the corresponding $\pi'$, equivalent to $\pi$, such that $i_a$ appears before $i_b$ in
$\pi'$ if there is an edge $(a,b)$ in $H_{\pi}(f)$, oriented from $a$ to $b$.
\end{proof}

\begin{example}
We illustrate this correspondence with the following example.  Let $G$ be a $4$-cycle
with vertices $1,\ldots ,4$, and let $\pi=(1,2,1,3)$.  Then $H_{\pi}$ has the four
vertices $1,11,2,3$, where $11$ represents the vertex corresponding to the second
$1$ in $\pi$.  There is an edge $3\rightarrow 1$, which becomes oriented toward $1$
in the acyclic orientation $\mathcal O_{\pi}$.
\end{example}

The next theorem provides an upper bound on the number of different systems 
$f^{\pi}:\kk^n\longrightarrow \kk^n$ one can obtain from composing the coordinate
functions of an $n$-tuple $f=(f^1,\ldots ,f^n)$, up to $t$ of them at a time.

\begin{theorem}\label{thm:gral}
Let $f=(f^1,\ldots ,f^n)$ be a system of local functions on $\kk^n$, and let
$F_{W_t}(f)=\{f^{\pi}|\pi\in W_t\}$.  Then
$$
|F_{W_t}(f)|\leq |\{{\rm Acyc}(H_{\pi}(f))|\pi\in W_t\}|
=\sum_{\pi\in W_t}|\{{\rm Acyc}(H_{\pi}(f)\}|.
$$
\end{theorem}

\begin{proof}
By Proposition \ref{main:lem1},
$ |W_t/\sim_{\overline{G}}| \leq |\{{\rm Acyc}(H_{\pi}(f)):\pi\in W_t\}|$.
But we have seen that if $\pi \sim_{\overline{G}} \pi'$ then $f^{\pi} =
f^{\pi'}$. Hence 
\[|F_{W_{t}}(f)|\leq |W_t/\sim_{\overline{G}}| 
\leq |\{{\rm Acyc}(H_{\pi}(f)):\pi\in W_t\}|.\]
\end{proof}  

This result shows in particular that if $\pi\sim_{\bar{G}}\pi'$, 
then the two systems $f^{\pi}$ and $f^{\pi'}$
are equal.
The following example shows that the upper bound in the theorem
is not attained in general.  

\begin{example}
Let $f = (x_{2}x_{3},x_{1}x_{3},0):\kk^{3}\to \kk^{3}$. Then $\Phi(f)$
does not contain any edges. Let $\pi = (3,2,1), \pi'=(3,1,2)$. Then
$\pi \nsim_{\bar{G}} \pi'$. However, 
$$
f^{\pi}=f^{1}\circ f^{2}\circ f^{3} =0=
f^{2}\circ f^{1}\circ f^{3} = f^{\pi'}.
$$
\end{example}

\begin{cor}
If $\pi \in S_{n}$, then $\{{\rm Acyc}(H_{\pi}(f)):\pi\in W_n\} = \{{\rm
Acyc}(\bar{G})\}$. Thus in this case we recover the upper bound for
the number of different SDS
obtained in \cite{bmr}.
\end{cor}

If we restrict ourselves to SDS, this bound is known to be sharp.  For
general systems this seems to be a substantially more difficult question.





\section{Finite Systems and SDS}

The approach to the study of finite systems taken in this paper was originated
in \cite{LP}, motivated by the desire to better understand sequential dynamical
systems.  

Recall that an SDS $f:\kk^n\longrightarrow \kk^n$ is given by a graph
$Y$ with $n$ vertices, functions $f^i:\kk^n\longrightarrow \kk^n$, which change only
the $i$th coordinate and take as input those coordinates connected to $i$
in the graph $Y$.  These functions are then composed according to an update
schedule given by a permutation $\pi\in S_n$.  That is, 
$$
f=f^{\pi(n)}\circ\cdots \circ f^{\pi(1)}.
$$
The functions $f^i$ are required to be symmetric in their inputs, that is,
permuting the inputs does not change the value of the function.  

In this paper we study $n$-tuples of functions $f^i:\kk^n\longrightarrow \kk^n$,
which change only the $i$th coordinate, without any further restrictions.
In particular, we do not suppose the a priori existence of a graph $Y$, that
governs the dependencies among these functions.  The Galois correspondence $\Psi$,
constructed in \cite{LP}, provides such a graph when needed.  And, as for SDS, it
is the invariants of this graph that determine many properties of the $n$-tuple and
finite systems derived from it.  It shows that even for SDS it is not necessary
to explicitly include the dependency graph $Y$ in the data defining an SDS.  The
Galois correspondence also shows that in general there will be more than one 
system whose dependency relations are modeled by a given graph.
There is, however, a natural
choice, which also provides a definition for the linearization of a system. 

An important
theoretical result, proved in \cite{R}, gives a sharp upper bound on the number
of different SDS that can be obtained by varying the update schedule over all of $S_n$.
The proof of this result assumes that all vertices of $Y$ that have the same degree
also have the same local function attached to them.
Theorem \ref{thm:gral} generalizes this upper bound by removing the restrictions on
the local functions $f^i$ and on the graph $Y$.  More importantly, it removes the
restriction that the update schedule be given by a permutation.  Thus, the upper bound
holds for compositions of the coordinate functions, which allows for repetitions of the
functions, and does not require that all functions are actually used.

These generalizations suggest that a more relaxed definition of SDS can still lead
to a class of systems about which one can prove theorems like the above upper bound.
Such a definition is proposed in \cite{LP2}, where a category of more general SDS
is developed.  

In \cite{R2} an upper bound for dynamically
non-equivalent SDS is given. In general, two maps $f,g: \kk^{n} \to
\kk^{n}$ are \emph{dynamically equivalent} if there exists a bijection
$\varphi: \kk^{n} \to \kk^{n}$ such that 
\[ g = \varphi \circ f \circ \varphi^{-1}.\]
This upper bound relies on the fact that conjugacy yields an SDS with
the same graph and local functions. This is not true for the  general
systems discussed in this paper as the following example shows.  Thus,
this upper bound holds exactly for the class of SDS.

\begin{example}
Let $f = (0,x_{3},x_{2}):\kk^{3}\to \kk^{3}$. Then $\Phi(f)$ is the
graph on three vertices $1,2,3$ with edges $(1,2),(1,3)$. Hence there
are only two functionally non-equivalent systems which correspond to
the permutations $id = (1 2 3)$ and $(3 2 1)$, that is, the systems
$f^{3}\circ f^{2}\circ f^{1}$ and $f^{1}\circ f^{2}\circ f^{3}$.  These
two systems have the state spaces in Figure \ref{states}.

\begin{figure}[!h]
\centering 
\includegraphics[totalheight=3cm]{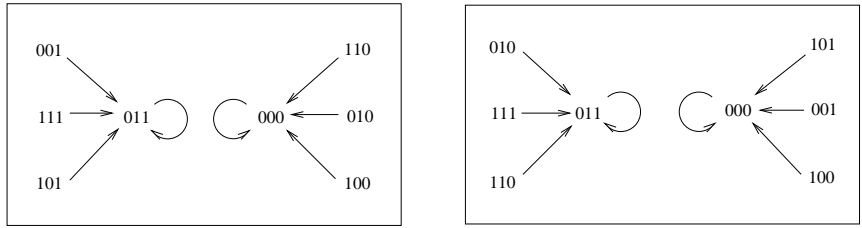}\label{fig:sys1}
\caption{The state spaces of $f^{3}\circ f^{2}\circ f^{1}$ and  
$f^{1}\circ f^{2}\circ f^{3}$.}   
\label{states}
\end{figure}

Nevertheless, if we let $\varphi = (2 1 3)$, then the system $\varphi\circ 
f^{id} \circ \varphi^{-1}$ has the state space in Figure \ref{conjugate}.

\begin{figure}[!htp]
\centering 
\includegraphics[totalheight=3cm]{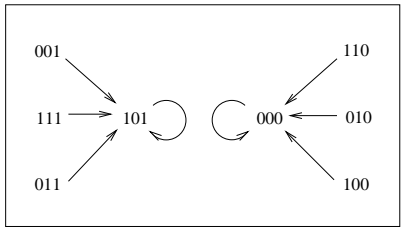}
\caption{The state space of $\varphi\circ f^{3}\circ f^{2}\circ
f^{1}\circ \varphi^{-1}$.}   
\label{conjugate}
\end{figure}

As expected, this state space is isomorphic to the state space of
$f^{id}$ but it is not equal to any of the two possible state spaces
given in Figure \ref{states}.  In fact, it is not equal to the state
space of any system obtained by composing the functions $f^1,f^2,f^3$
according to any word.  This is easily seen because, for example,
the state $(1,0,1)$ is sent to itself, but $f^1$ is the zero function.
So $f^1$ cannot be involved.  On the other hand, the first coordinate
of other states changes, so there must be a function involved that 
changes the first coordinate.
\end{example}

\vspace{0.1cm}
\section{Stably isomorphic systems}

In this section we answer a question that suggests itself naturally
from the Galois correspondence.  Given two systems that have the
same dependency graphs, are they stably isomorphic, that is, do they
have the same limit cycle structure in their state spaces?  Recall
that the state space of a system is a directed graph, whose vertices
are the states of the system, that is, all binary strings of a given
length.  Directed edges correspond to system transitions.

\begin{definition}
Let $f_{1}, f_{2}$ be systems with state spaces $\mathcal{S}_{f_{i}}$
and subdigraphs of limit cycles $\mathcal{L}_{f_{i}}$. We call $f_{1}$
and $f_{2}$ \emph{stably isomorphic} if there exists a digraph
isomorphism between $\mathcal{L}_{f_{1}}$ and $\mathcal{L}_{f_{2}}$.
\end{definition}

Let $f = (f^1,\dots,f^n)$ and $g = (g^1,\dots,g^n)$ 
be two systems with the same dependency graphs, that is, $\Phi(f) = \Phi(g)$.  
One can now 
ask if $f^{\pi(n)}\circ \cdots \circ f^{\pi(1)}$ is
stably isomorphic to 
$g^{\pi(n)}\circ \cdots \circ g^{\pi(1)}$
for all $\pi\in W_t$. 
In this section we show that this is not true, 
by providing a counterexample.

\begin{example}
Let $f = (f^1,f^2,f^3)$ and $g = (g^1,g^2,g^3)$ be two triples of
functions with
\begin{eqnarray*}
f^1(x_1,x_2,x_3) &=& (1+x_2x_3,x_2,x_3), \\
f^2(x_1,x_2,x_3) &=& (x_1,x_1,x_3), \\
f^3(x_1,x_2,x_3) &=& (x_1,x_2,1+x_1), \\
g^1(x_1,x_2,x_3) &=& (x_1+x_2+x_3,x_2,x_3), \\
g^2(x_1,x_2,x_3) &=& (x_1,1+x_1+x_2,x_3), \\
g^3(x_1,x_2,x_3) &=& (x_1,x_2,x_1+x_3). 
\end{eqnarray*}
By using Theorem \ref{main:thm1}, 
it is easy to see that $\Phi(f) = \Phi(g)$ is the graph on vertices
$1,2,3$, with the single edge $(2,3)$.

Let $f = f^3 \circ f^2 \circ f^1$ and $g = g^3 \circ g^2 \circ g^1$. 
The state spaces of $f$ and $g$ are given in Figure \ref{state}.
 
\begin{figure}[!htp] 
\centering 
\includegraphics[totalheight=3cm]{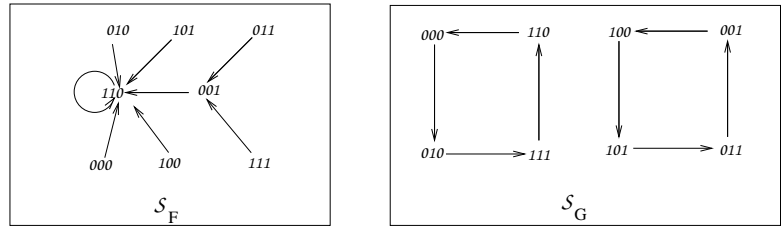} 
\caption{The state spaces of $f$ and $g$.} \label{state} 
\end{figure}

{}From Figure \ref{state}, it is clear that $f$ and $g$ are not 
stably isomorphic. 
\end{example}


\end{document}